\documentclass[11pt]{article}
\usepackage[french, english]{babel}
\usepackage{amsmath,amssymb,amsthm,mathrsfs,url}
\usepackage{a4wide}
\usepackage{xypic}
\newtheorem{thm}{Theorem}
\newtheorem{theoreme}[thm]{Th\'eor\`eme}
\newtheorem{prop}[thm]{Proposition}
\newtheorem{cor}[thm]{Corollary}
\newtheorem{corollaire}[thm]{Corollaire}
\newtheorem{conj}[thm]{Conjecture}
\newtheorem{rem}[thm]{Remark}

\newcommand{\field}[1]{\mathbb{#1}}
\newcommand{\Q}{\field{Q}}

\newcommand{\F}{\field{F}}

\newcommand{\T}{\field{T}}

\title{A non-solvable Galois extension of $\Q$ ramified at 2 only}
\author{Lassina Demb\'el\'e}\date{}

\begin{document}
\maketitle
\begin{quote}\small\it \`A la m\'emoire de ma s{\oe}ur jumelle Fatouma. D\'ej\`a vingt ans que tu es partie \end{quote}

\abstract{In this paper, we show the existence of a non-solvable Galois extension of $\Q$ which is unramified outside $2$. The extension $K$ we construct has degree $2251731094732800=2^{19}(3\cdot 5\cdot 17\cdot 257)^2$ and has root discriminant $\delta_K <2^{\frac{47}{8}}=58.68...$, and is totally complex.}

{\selectlanguage{french}
\abstract{Dans cet article, nous d\'emontrons l'existence d'une extension galoisienne non r\'eso\-luble de $\Q$ ramifi\'ee seulement en $2$. L'extension $K$ que nous construisons est de degr\'e $2251731094732800=2^{19}(3\cdot 5\cdot 17\cdot 257)^2$ et de discriminant normalis\'e $\delta_K < 2^{\frac{47}{8}}=58,68...$, et est totalement complexe.}
}

\section*{\bf Version fran\c{c}aise abr\'eg\'ee}
La conjecture suivante est propos\'ee dans Gross~\cite{gross1}.
\begin{conj}\label{conj:nonramifiee} Pour tout nombre premier $p$, il existe une extension galoisienne non r\'esoluble de $\Q$ ramifi\'ee seulement en $p$.
\end{conj}
Ce r\'esultat est connu lorsque $p \ge 11$.  En effet, Serre~\cite{serre3} montre que pour un tel nombre premier $p$, on trouve $k=12,16,18,20,22$ ou $26$ tel que la repr\'esentation galoisienne r\'esiduelle $\bar{\rho}_{k,p}$ mod $p$ attach\'ee \`a l'unique forme parabolique de niveau 1 et de poids $k$, \`a coefficients entiers, est absolument irr\'eductible.  Par~\cite[Chap. IV]{serre5}, le corps fixe de $\ker\bar{\rho}_{k, p}$ est donc une extension non r\'esoluble de $\Q$ qui est non ramifi\'ee en dehors de $p$.

Dans cet article, nous \'etablissons cette conjecture pour $p=2$. L'extension est construite \`a partir de repr\'esentations galoisiennes attach\'ees aux formes modulaires de Hilbert de niveau 1, de poids parall\`ele 2 et \`a coefficients dans $\overline{\F}_2$, sur le sous-corps totalement r\'eel maximal $F$ du corps cyclotomique $\Q(\zeta_{32})$. Nous d\'emontrons le th\'eor\`eme suivant:

\begin{theoreme}\label{thm:nonsoluble} Il existe deux $\mathbf{SL}_2(\F_{2^8})$-extensions $E$ et $E'$ de $F$ ramifi\'ees en l'unique id\'eal premier divisant $2$. Les extensions $E$ et $E'$ sont galoisiennes sur $\Q(\sqrt{2})$, avec groupe de Galois $\mathbf{SL}_2(\F_{2^8})\cdot 4$, et sont interchang\'ees par $\mathrm{Gal}(\Q(\sqrt{2})/\Q)$.
\end{theoreme}

\noindent
On en d\'eduit:

\begin{corollaire} Il existe une extension galoisienne non r\'esoluble $K$ de $\Q$ qui est ramifi\'ee seulement en $2$, de groupe de Galois $\mathbf{SL}_2(\F_{2^8})^2\cdot 8$.
\end{corollaire}

Une \'etude locale des repr\'esentations galoisiennes \`a partir desquelles l'extension $K$ \`a \'et\'e construite nous permet de borner son discriminant. On obtient ainsi le r\'esultat suivant. (Le ``discriminant normalis\'e'' d'une extension finie  $E$  de  $\Q$  est $|d_E|^{1/[E:\Q]}$, o\`u  $d_E$  est le discriminant de  $E$.)

\begin{prop}\label{borne_disc} Le discriminant  normalis\'e $\delta_K$ de l'extension $K$ est $<2^{\frac{47}{8}}=58,68...$
\end{prop}
La Proposition~\ref{borne_disc} implique que l'extension $K$ ne peut \^etre totalement r\'eelle; sinon, on aurait $\delta_K>60, 83...$, la borne inf\'erieure d'Odlyzko pour un corps totalement r\'eel de ce degr\'e, que nous avons \'evalu\'ee par les formules de Poitou~\cite{diaz1}. Elle est donc totalement complexe, \'etant donn\'e qu'elle est galoisienne.

\begin{rem}\rm En fait, la borne de la Proposition~\ref{borne_disc} peut \^etre abaiss\'ee \`a $\delta_K\le 55,39....$ \`A ce sujet, nous r\'ef\'erons au compl\'ement de Jean-Pierre Serre qui suit.
\end{rem}

\setcounter{thm}{0}
\section{\bf Introduction}
In this paper, we prove that there exists a non-solvable finite Galois extension of $\Q$ which is ramified at 2 only. We construct this extension by using the Galois representations attached to Hilbert modular forms over the maximal totally real subfield of the cyclotomic field $\Q(\zeta_{32})$. This settles the following conjecture, proposed in Gross~\cite{gross1}, for the prime $p=2$.

\begin{conj}\label{conj:unramified} For any prime number $p$, there is a finite non-solvable Galois extension $K$ of $\Q$ ramified at $p$ only.
\end{conj}

For primes $p\ge 11$, one knows how to construct extensions satisfying Conjecture~\ref{conj:unramified}. Indeed, Serre~\cite{serre3} shows that for such a prime $p$, there is $k=12,16,18,20,22$ or $26$ such that the residual Galois representation $\bar{\rho}_{k,p}$ mod $p$ associated to the unique cuspidal form of level 1 and weight $k$, with integral coefficients, is absolutely irreducible. By~\cite[Chap. IV]{serre5}, the fixed field of $\ker\bar{\rho}_{k, p}$ is then a non solvable extension of $\Q$ unramified away from $p$.

As for primes $\le 7$, the first case of the Serre conjecture~\cite{serre1} was proved and later published by Tate~\cite{tate1}, for $p=2$, by simply ruling out the existence of mod 2 irreducible representations of the absolute Galois group $\mathrm{Gal}(\bar{\Q}/\Q)$, unramified away from 2. His results were later extended to the primes 3 and 5 by \cite{serre4} and \cite{brueggeman1} respectively, assuming GRH in the latter case. By Khare and Wintenberger~\cite{khare1}, this is now true unconditionally for  $p=5, 7$.

However, the general belief is that one can still solve Conjecture~\ref{conj:unramified} by working with automorphic forms over algebraic groups of higher rank; for instance, exceptional groups. Unfortunately, it is not yet clear how to attach Galois representations to such automorphic forms in general. Underscoring this, Gross~\cite{gross2} developed the theory of algebraic modular forms and conjectured the existence of Galois representations attached to irreducible Hecke constituents. Computational results for the exceptional group $G_2$ over $\Q$ based on his conjectures, led Lansky and Pollack~\cite{lansky1} to predict the existence of a $G_2(\F_5)$-extension of $\Q$ that is ramified at 5 only, thus providing further evidence for Conjecture~\ref{conj:unramified}. 

Our approach, which was suggested to us by Gross, relies on fixing the underlying group and enlarging the base field instead; as for us, it was much easier to study groups of higher rank that are not absolutely simple. The extension $K$ we construct has degree $2^{19}(3\cdot 5\cdot 17\cdot 257)^2\sim 2\times 10^{15}$ and has root discriminant $\delta_K <2^{\frac{47}{8}}=58.68..$. Thus it is totally complex, being Galois over $\Q$. To the best of our knowledge, this is the largest known totally complex field with such a low root discriminant. It would be interesting to know whether this is an isolated case or if there are infinite towers of totally complex fields with their minimal root discriminant in this magnitude. Indeed, the current upper bound for such towers is 82.2, which was obtained by Hajir and Maire~\cite{hajir1, hajir2}. With the Hilbert Modular Forms Package being currently implemented in Magma~\cite{magma1}, we hope to settle the remaining cases of Conjecture~\ref{conj:unramified} in the near future.

\bigskip
\noindent{\bf Acknowledgements.} I would like to thank Dick Gross for suggesting this question, and for his extreme generosity, enthusiasm and encouragement. I would like to thank Jean-Pierre Serre for carefully reading an earlier version of this note, and for making numerous suggestions that help improve the presentation. I would like to thank the Magma group at the University of Sydney for their support, especially the assistance of Steve Donnelly. I would also like to thank Kevin Buzzard, Fred Diamond and David Roberts for helpful email exchanges, as well as Gabor Wiese for useful conversations. This project was funded by a grant of SBF/TR 45 of the Deutsche Forschungsgemeinschaft.

\section{\bf A non-solvable extension of $\Q$ ramified at 2 only}\label{sec1}  
Let $F$ be the maximal totally real subfield of the cyclotomic field $\Q(\zeta_{32})$, and $\mathcal{O}_F$ its ring of integers. It is generated by the element $\beta:=\zeta_{32}+\zeta_{32}^{-1}$ with minimal polynomial $x^8 - 8x^6 + 20x^4 - 16x^2 + 2$. We fix the integral basis $1,\beta,\cdots,\beta^7$ of $F$, and we let $\sigma$ be the cyclic generator of $\mathrm{Gal}(F/\Q)$ given by $(\beta\mapsto -\beta^3 + 3\beta)$. We let $\alpha$ be a cyclic generator of $\F_{2^8}^\times$, the unit group in $\F_{2^8}$.

\begin{thm}\label{thm:nonsolvable} There exist two $\mathbf{SL}_2(\F_{2^8})$-extensions $E$ and $E'$ of $F$ ramified at the unique prime ideal above $2$ only. The extensions $E$ and $E'$ are  both Galois over $\Q(\sqrt{2})$, with Galois group $\mathbf{SL}_2(\F_{2^8})\cdot 4$, and are interchanged by $\mathrm{Gal}(\Q(\sqrt{2})/\Q)$.
\end{thm}

\begin{proof} Let $S_2(1,\F_2)$ be the space of $\mod 2$ Hilbert cusp forms of level $1$ and parallel weight 2 over $F$. Let $\T$ be the Hecke algebra over $\F_2$ generated by the operator $T(\mathfrak{p})$, where $\mathfrak{p}$ runs over all the primes in $F$. We computed the space $S_2(1,\F_2)$ and the action of $\T$ on it using the Hilbert Modular Forms Package in Magma~\cite{magma1}. It has two nonzero irreducible Hecke constituents which are both 8 dimensional over $\F_2$. The action of $\T$ is completely determined by the operators $T(\mathfrak{p}_{31}^i)$ at the primes above $31$, which splits completely in $F$ into 8 distinct primes $\mathfrak{p}_{31}^i$, $i=1,\,\ldots,\,8$ (see Table~\ref{table1} for notations). The common characteristic polynomial of those operators together with the one of $T(\mathfrak{p}_2)$, where $\mathfrak{p}_2$ is the unique prime above 2, are given by 
\begin{eqnarray*}
\mathrm{charpoly}\left(T(\mathfrak{p}_{2})\right)&=&x^{41}(x^2+x+1)^8\mod 2\\
\mathrm{charpoly}\left(T(\mathfrak{p}_{31}^1)\right)&=&x^{41}(x^8 + x^4 + x^3 + x + 1)(x^8 + x^6 + x^5 + x^2 + 1)\mod 2.
\end{eqnarray*}
Let $f$ and $f'$ be the newforms whose first few Hecke eigenvalues are listed in Table~\ref{table1}. Their $\mathrm{Gal}(\F_{2^8}/\F_2)$-conjugacy classes determine the two nonzero constituents of $S_2(1,\F_2)$. We recall that to give a newform $f\in S_2(1,\overline{\F}_2)$ is equivalent to giving a maximal ideal $\mathfrak{m}_f\subset\T$; and that the association ($f\mapsto\mathfrak{m}_f$) is a bijection between $\mathrm{Gal}(\overline{\F}_2/\F_2)$-conjugacy classes of newforms and maximal ideals in $\T$.

\medskip
Let $\mathfrak{m}_f,\,\mathfrak{m}_{f'}\subset \T$ be the maximal ideals associated to $f$ and $f'$ respectively, and let $\theta_f:\,\T\to\T/\mathfrak{m}_f=\F_{2^8}$ and $\theta_{f'}:\,\T\to \T/\mathfrak{m}_{f'}=\F_{2^8}$ be the corresponding ring homomorphisms.  By work of Rogawski-Tunnell, Ohta and Carayol~\cite{rogawski1, ohta1, carayol1}, completed by Taylor and Jarvis~\cite{taylor1, jarvis1, taylor2}, there are Galois representations
$$\bar{\rho}_f,\,\bar{\rho}_{f'}:\,\mathrm{Gal}(\overline{F}/F)\to\mathbf{SL}_2(\F_{2^8})$$ such that
$\mathrm{Tr}\left(\bar{\rho}_{f}(\mathrm{Frob}_{\mathfrak{p}})\right)=\theta_f(T(\mathfrak{p}))$ and $\mathrm{Tr}\left(\bar{\rho}_{f'}(\mathrm{Frob}_{\mathfrak{p}})\right)=\theta_{f'}(T(\mathfrak{p}))$, for all prime $\mathfrak{p}$.  From the orders of Frobenii provided in Table~\ref{table1}, we see that $\bar{\rho}_{f}$ and $\bar{\rho}_{f'}$ are surjective. This proves the first part of Theorem~\ref{thm:nonsolvable} with the two extensions $E$ and $E'$ being the fixed fields of $\ker(\bar{\rho}_{f})$ and $\ker(\bar{\rho}_{f'})$ respectively. (We recall that the extensions $E$ and $E'$ only depend on the $\mathrm{Gal}(\F_{2^8}/\F_2)$-conjugacy classes of $f$ and $f'$, or equivalently $\mathfrak{m}_f$ and $\mathfrak{m}_{f'}$, respectively.)

\medskip
 Let $\tau$ be the cyclic generator of $\mathrm{Gal}(\F_{2^8}/\F_2)$ given by $(\tau:\F_{2^8}\to\F_{2^8},\,\alpha\mapsto\alpha^2)$. It is not hard to see that, for the primes listed in Table~\ref{table1}, 
$$a_{\sigma^2(\mathfrak{p})}(f)=\tau^2(a_{\mathfrak{p}}(f))\,\mbox{\rm and }\, a_{\sigma^2(\mathfrak{p})}(f')=\tau^2(a_{\mathfrak{p}}(f')).$$ And since those primes determine $f$ and $f'$, these identities extend to all primes $\mathfrak{p}$. This means that the action of $\mathrm{Gal}(F/\Q(\sqrt{2}))=\langle\sigma^2\rangle$ preserves the $\mathrm{Gal}(\F_{2^8}/\F_2)$-conjugacy classes of $f$ and $f'$. Or equivalently, that
$\sigma^2(\mathfrak{m}_f)=\mathfrak{m}_{\tau^2(f)}=\mathfrak{m}_f$ and $\sigma^2(\mathfrak{m}_{f'})=\mathfrak{m}_{f'}$. From this, we conclude that $E$ and $E'$ are Galois over $\Q(\sqrt{2})$ with the same Galois group $\mathbf{SL}_2(\F_{2^8})\cdot 4$. This proves the second part of Theorem~\ref{thm:nonsolvable}.

\medskip
Finally, we observe that  $a_{\sigma(\mathfrak{p})}(f)=a_{\mathfrak{p}}(f')$, for any prime $\mathfrak{p}$, which implies that $\sigma(\mathfrak{m}_f)=\mathfrak{m}_{f'}$. Therefore, the action of $\mathrm{Gal}(F/\Q)$ permutes the $\mathrm{Gal}(\F_{2^8}/\F_2)$-conjugacy classes of $f$ and $f'$. Combining this with the observation above, we see that $\mathrm{Gal}(\Q(\sqrt{2})/\Q)$ interchanges $E$ and $E'$, which concludes the proof of Theorem~\ref{thm:nonsolvable}.
\end{proof}

\begin{table}
\small
\begin{eqnarray*}
\begin{array}{|c|cccccccc|}\hline
\mathfrak{p}&\mathfrak{p}_{31}^{1}&\mathfrak{p}_{31}^{2}&\mathfrak{p}_{31}^{3}&\mathfrak{p}_{31}^{4}&\mathfrak{p}_{31}^{5}&\mathfrak{p}_{31}^{6}&\mathfrak{p}_{31}^{7}&\mathfrak{p}_{31}^{8}\\\hline\hline
a_\mathfrak{p}(f)&\alpha^{100}& \alpha^{19}& \alpha^{145}& \alpha^{76}& \alpha^{70}& \alpha^{49}& \alpha^{25}& \alpha^{196}\\
\mathrm{ord}(\bar{\rho}_{f}(\mathrm{Frob}_\mathfrak{p}))&257& 255& 257& 255& 257& 255& 257& 255\\\hline

a_\mathfrak{p}(f')&\alpha^{196}& \alpha^{100}& \alpha^{19}& \alpha^{145}& \alpha^{76}& \alpha^{70}& \alpha^{49}& \alpha^{25}\\
\mathrm{ord}(\bar{\rho}_{f'}(\mathrm{Frob}_\mathfrak{p}))&  255& 257& 255& 257& 255& 257& 255& 257\\\hline
\multicolumn{9}{}{}\\\hline

\mathfrak{p}&\mathfrak{p}_{97}^{1}&\mathfrak{p}_{97}^{2}&\mathfrak{p}_{97}^{3}&\mathfrak{p}_{97}^{4}&\mathfrak{p}_{97}^{5}&\mathfrak{p}_{97}^{6}&\mathfrak{p}_{97}^{7}&\mathfrak{p}_{97}^{8}\\\hline\hline
a_\mathfrak{p}(f)&\alpha^{23}& \alpha& \alpha^{92}& \alpha^4& \alpha^{113}& \alpha^{16}& \alpha^{197}& \alpha^{64}\\
\mathrm{ord}(\bar{\rho}_{f}(\mathrm{Frob}_\mathfrak{p}))&257& 51& 257& 51& 257& 51& 257& 51\\\hline
a_\mathfrak{p}(f')& \alpha^{64}& \alpha^{23}& \alpha& \alpha^{92}& \alpha^4& \alpha^{113}& \alpha^{16}& \alpha^{197}\\
\mathrm{ord}(\bar{\rho}_{f'}(\mathrm{Frob}_\mathfrak{p}))&51&257& 51& 257& 51& 257& 51& 257\\\hline
\end{array}
\end{eqnarray*}
\begin{eqnarray*}
\mathfrak{p}_{31}&:=&([1, 2, 0, -4, 0, 1, 0, 0]),\\ 
\mathfrak{p}_{31}^i&:=&\sigma^{i-1}(\mathfrak{p}_{31}),\,\,i=1,\ldots,8.\\
\mathfrak{p}_{97}&:=&([1, -12, -4, 19, 1, -8, 0, 1]),\\ 
\mathfrak{p}_{97}^i&:=&\sigma^{i-1}(\mathfrak{p}_{97}),\,\,i=1,\ldots,8.\\
\end{eqnarray*}
\caption{\small Mod $2$ Hilbert newforms of weight $2$ and level $1$ over $F=\Q(\zeta_{32})^+$.} 
\label{table1}
\end{table}

\begin{cor} There exists a finite non-solvable Galois extension $K$ of $\Q$ that is ramified at $2$ only, with Galois group $\mathbf{SL}_2(\F_{2^8})^2\cdot 8$.
\end{cor}

\begin{proof} Let $K$ be the Galois closure of $E$ over $\Q$. By Theorem~\ref{thm:nonsolvable}, $K$ is the compositum of $E$ and $E'$, and only ramifies at $2$ by construction. So it only remains to find its Galois group.

By Galois theory and the fact that $E\cap E'$ is Galois over $F$, we have
\begin{eqnarray*}
\mathrm{Gal}(K/E)&=&\mathrm{Gal}(EE'/E)\cong\mathrm{Gal}(E/E\cap E')=\mathrm{Gal}(E'/E\cap E')\\
&\lhd&\mathrm{Gal}(E/F)=\mathbf{SL}_2(\F_{2^8}).
\end{eqnarray*} Since $E$ is not Galois over $\Q$ and the only normal subgroups of $\mathbf{SL}_2(\F_{2^8})$ are $1$ and itself, we must have $E\cap E'=F$ and $\mathrm{Gal}(K/E)=\mathrm{Gal}(K/E')=\mathbf{SL}_2(\F_{2^8})$. Thus, the fields $E$ and $E'$ are disjoint over $F$ and $\mathrm{Gal}(K/F)=\mathbf{SL}_2(\F_{2^8})^2$, which implies that
\begin{eqnarray*}
\mathrm{Gal}(K/\Q)&\cong&\mathrm{Gal}(K/F)^2\cdot 8=\mathbf{SL}_2(\F_{2^8})^2\cdot 8.
\end{eqnarray*}
\end{proof}

\section{\bf An estimate for the root discriminant} 

In this section, we use the modularity of the Galois representations from which our extension $K$ of Section~\ref{sec1} was constructed in order to obtain an estimate for its root discriminant.

\begin{prop}\label{disc_bound} The root discriminant $\delta_K$ of the extension $K$ is $<2^{\frac{47}{8}}=58.68.. .$
\end{prop}

\begin{proof} Let $\mathfrak{P}_1,\ldots,\mathfrak{P}_{g}$ be the primes in $\mathcal{O}_K$ above $\mathfrak{p}_2$ so that
$$\mathfrak{p}_2=\prod_{i=1}^{g}\mathfrak{P}_i^{e_i}.$$ Since $K$ is Galois over $F$, the group $\mathrm{Gal}(K/F)$ acts transitively on the set of those primes. We let  $e$ and $s$ be the common ramification index and residue field degree, respectively, so that  $e_i=e$ and $esg=[K:F]$. Let $\mathfrak{P}$ be any of the primes above $\mathfrak{p}_2$, and $\mathfrak{p}=\mathcal{O}_E\cap\mathfrak{P}$ and  $\mathfrak{p}'=\mathcal{O}_{E'}\cap\mathfrak{P}$. Let $E_\mathfrak{p}$, $E_{\mathfrak{p}'}'$ and  $K_\mathfrak{P}$ be the completions of $E$, $E'$ and  $K$ at $\mathfrak{p}$, $\mathfrak{p}'$ and $\mathfrak{P}$ respectively. Then $K_{\mathfrak{P}}$ is the compositum of $E_{\mathfrak{p}}$ and $E_{\mathfrak{p}'}'$. From the characteristic polynomial of $T(\mathfrak{p}_2)$, we see that the form $f$ is ordinary at $\mathfrak{p}_2$, with $a_{\mathfrak{p}_2}(f)$ being a generator of $\F_4^\times$. By Wiles~\cite[Theorem 2]{wiles1} it follows that the restriction of $\bar{\rho}_f$ to the decomposition group at $\mathfrak{p}_2$ is of the form
$$\bar{\rho}_f|D_{\mathfrak{p}_2}\sim\begin{pmatrix}\chi&*\\ 0&\chi^{-1}\end{pmatrix},$$ where $\chi$ is an unramified character of order 3. It also follows that 
$$\bar{\rho}_f|I_{\mathfrak{p}_2}\sim\begin{pmatrix}1&*\\ 0&1\end{pmatrix}.$$ From this and the fact that the extensions $E_{\mathfrak{p}}$ and $E_{\mathfrak{p}'}'$ are {\it peu ramifi\'ees} in the sense of Serre~\cite[sec. 2]{serre1} (see also~\cite{breuil1} for a group scheme theoretic definition), it follows that  $K_{\mathfrak{P}}=L(\sqrt{x_1},\ldots,\sqrt{x_{m}})$, where $e=2^m$, $L$ is the unique unramified extension of degree 3 of $F_{\mathfrak{p}_2}$ contained in $K_{\mathfrak{P}}$ and $x_i\in\mathcal{O}_L^\times/\left(\mathcal{O}_L^\times\right)^2$. And so, the Galois group $\mathrm{Gal}(K_{\mathfrak{P}}/L)$ has $2^m-1$ quadratic characters, whose conductors divide $\mathfrak{p}_L^{16}$. Therefore, by the discriminant-conductor formula~\cite[Chap. VI]{serre2}, we get that the local discriminant $d_{K_\mathfrak{P}/L}$ divides $\mathfrak{p}_L^{16(2^m-1)}.$ Equivalently, this means that
$d_{K_\mathfrak{P}/F_{\mathfrak{p}_2}}$divides $\hat{\mathfrak{p}}_2^{16(2^m-1)},$ where $\hat{\mathfrak{p}}_2$ is the maximal ideal in $F_{\mathfrak{p}_2}$. Taking the product over all primes then yields that the global discriminant $d_{K/F}$ divides $\mathfrak{p}_2^{16gs(2^m-1)}=\mathfrak{p}_2^{2[K:\Q](1-1/2^m)}.$ From the relation
\begin{eqnarray*}
d_{K}=d_{F}^{[K:F]}\mathrm{N}_{F/\Q}(d_{K/F}),
\end{eqnarray*} it then follows that $d_{K}$ divides $2^{31[K:F]}\times 2^{2[K:\Q](1-1/2^m)}$, and hence
\begin{eqnarray*}
\delta_{K}\le\delta_{F}2^{2(1-\frac{1}{2^m})}=2^{\frac{31}{8}}2^{2(1-\frac{1}{2^m})}=2^{\frac{47}{8}-\frac{1}{2^{m-1}}}<2^{\frac{47}{8}}.
\end{eqnarray*}
\end{proof}
From Proposition~\ref{disc_bound}, we see that $K$ cannot be totally real; otherwise, we would have $\delta_K>60.83$, the Odlyzko bound for a totally real field of this degree, estimated using Poitou's formulas in~\cite{diaz1}. Therefore it must be totally complex, being Galois over $\Q$. 

\begin{rem}\rm The bound in Proposition~\ref{disc_bound} can be lowered to $\delta_K\le 55.39...$ To this end, we refer to the supplement written by Jean-Pierre Serre.
\end{rem}

\end{document}


\begin{center}
 Un complément à la Note de L.Dembélé
 
  ``A non-solvable Galois extension of  \;\Q \; 
 ramified at 2 only''
  \end{center}

  \begin{center}
  Jean-Pierre Serre, Collège de France, 3 rue d'Ulm, 75005 Paris
  \end{center}
  
  La Note de L.Dembélé décrit une extension galoisienne  $K/\Q$
  dont le degré  $n$  est grand ($n = 2^{19}.3^2.5^2.17^2.257^2 = 2251731094732800$), et qui a plusieurs propriétés intéressantes:
  
   1) Elle est non ramifiée en dehors de 2.
   
    2) Son groupe de Galois n'est pas résoluble: il contient comme sous-groupe d'indice 8 le produit direct de deux exemplaires de  $\SL_2(\F_{256}).$
    
    3) Son discriminant normalisé $\delta_K = |d_K|^{1/n}$ est petit: il est majoré
    par 58,688..., alors que les meilleurs exemples connus, dus à F.Hajir et C.Maire, étaient voisins de 82.
    
     En fait, lorsque l'on examine plus en détail la ramification du corps $K$, on obtient
     un résultat encore meilleur, à savoir $\delta_K \leqslant 55,394388...$. Voici
     pourquoi:
    
    \smallskip
        Je conserve les notations de la Note. Le point essentiel consiste à estimer, aussi précisément que possible, les conducteurs des caractères associés
        à l'extension multiquadratique  $K_{\mathfrak{P}}$ du corps  $L$; cette extension est  de
    degré  $2^{2m}$, avec $m \leqslant 8$.
    
       Soit $X $ le groupe des caractères en question, i.e. le dual de $ \Gal(K_{\mathfrak{P}}/L)$; c'est un sous-groupe  de  $Y = L^\times/L^{\times2}$. Le groupe  $Y$  est un $\F_2$-espace vectoriel de dimension 3.8 + 1 = 25. Il est muni d'une action du groupe $C$ d'ordre 3 qui est le groupe de Galois
    de l'extension $L/F_{\p_2}$. Cette action le décompose en deux morceaux:
    $$  Y = Y^+ \oplus Y^- ,$$
  où  $Y^+$  est la partie fixe, et $Y^-$ est l'unique supplémentaire stable. Si  $c$
  est un générateur de  $C$, $Y^-$ est l'ensemble des $y \in Y$  tels que $(1+c+c^2)y = 0.$ On peut donc voir $Y^-$ comme un $\F_4$-espace vectoriel. On vérifie facilement que  $\dim Y^+ = 9$ (car c'est le ``$Y$ '' du corps  $F_{\p_2}$), et que 
  la
  dimension de $Y^-$ sur $\F_4$ est égale à 8, de sorte que  $Y^-$  est d'ordre
  $4^8 = 2^{16}$.
  \smallskip
  
    L'intérêt de cette décomposition de  $Y$  en deux morceaux est que {\it le groupe $X$ associé à l'extension $K_{\mathfrak{P}}/L$ est contenu dans $Y^-$.} C'est évident lorsque l'on regarde
    le produit semi-direct  $\Gal(K_{\mathfrak{P}}/L).C$  comme plongé dans le produit  de deux copies de $\SL_2(\F_{2^8})$. [Noter que ceci entraîne que l'extension  $K_{\mathfrak{P}}/L$  est peu ramifiée,
    car les conducteurs des éléments de  $Y^-$ sont au plus égaux à $2e = 16. ]$
      Il faut ensuite regarder la filtration de  $X$  définie par le conducteur. Celle de
      $Y^-$  est donnée par un drapeau complet 
       $$ Y^-= Y_0 \supset Y_1 \supset  ... \supset Y_8 = 0, $$
    où la dimension sur $\F_4$ de chaque quotient successif $Y_i/Y_{i+1}$ est 1,
    et le conducteur de tout élément de  $Y_i - Y_{i+1}$ est $16-2i$. Lorsqu'on intersecte
    ceci avec  $X$, on obtient une filtration de $X$ avec quotient successifs, soit
    nuls, soit de dimension   1  sur  $\F_4$ . Si l'on appelle  $c_1,...,c_m$ les sauts de
    cette filtration, ordonnés par  $16 \geqslant c_1 > c_2 > ... > c_m > 0,$ on en conclut que $X$ contient:
    
      $4^m - 4^{m-1}$ éléments de conducteur $c_1,$
      
      $4^{m-1} - 4^{m-2}$ éléments de conducteur $c_2,$
      
      ...
      
      $ 4 - 1$ éléments de conducteur $c_m.$
      
      D'où la valuation du discriminant normalisé de $K_{\mathfrak{P}}/L$, à savoir
      \smallskip
      
   $v_L$(disc.norm) = $4^{-m}((4^m - 4^{m-1})c_1 + (4^{m-1} - 4^{m-2})c_2 + ... + (4 - 1)c_m)$,
    
    \n ce que l'on peut écrire sous la forme
    $v_L$(disc.norm) = $ c_1 - \epsilon, $ avec 
    $$ \epsilon = (c_1 - c_2)/4 + (c_2 - c_3)/4^2 + ... + c_m/4^m.$$
    Les $c_i - c_{i+1}$ sont pairs et $> 0.$  On en déduit
    $$ \epsilon \geqslant 2/4 + 2/4^2 + ... + 2/4^m = \frac{2}{3}(1 - 4^{-m}),$$
    d'où  \; $v_L$(disc.norm) $\leqslant c_1 - \frac{2}{3}(1 - 4^{-m}).$
      Sous cette forme, il est facile de prouver que la valeur maximale de $v_L$(disc.red)
      est atteinte pour  $c_1 = 16$ et $m=8$  et l'on obtient alors
      
     $v_L$(disc.norm) $ \leqslant 16 - \epsilon$  \; avec  $\epsilon = \frac{2}{3}(1 - 4^{-8}) = 21845/2^{15}.$
      Lorsque l'on utilise cette valeur pour estimer le discriminant normalisé du corps  $K$
      (sur $\Q$, cette fois), on obtient
        $\delta_K = 2^z,$
        avec 
        $$ z = 31/8 + 2 - \epsilon/8 = 1518251/262144 =  5,79166793...,$$
        d'où $\delta_K \leq 55,394388...$.